\documentclass{article}

\usepackage{amsthm}
\usepackage{amssymb}
\usepackage{amsxtra}
\usepackage{setspace}

\newtheorem{thm}{Theorem}[section]
\newtheorem{lem}{Lemma}[section]

\newtheorem{conj}{Conjecture}[section]
\newtheorem{question}{Question}

\title{Finding Almost Squares}
\author{Tsz Ho Chan}

\begin{document}
\maketitle
\begin{abstract}
We study short intervals which contain an ``almost square'', an integer $n$ that
can be factored as $n = ab$ with $a$, $b$ close to $\sqrt{n}$. This is related to the problem on distribution of $n^2 \alpha \pmod 1$ and the problem on gaps between sums of two squares.
\end{abstract}
\section{Introduction}

For any $a \geq 0$, there is a perfect square in the interval
$[a^2, (a+1)^2]$, namely $([a]+1)^2$ where $[x]$ stands for the
greatest integer smaller than $x$. Thus, for $x \geq 0$, the
interval $[x, x + 2\sqrt{x} + 1]$ always contains a perfect
square. Instead of perfect squares, we can look for integer $n$
that can be factored as $n = ab$ with $a$, $b$ close to
$\sqrt{n}$. We formulate the question as follow:
\begin{question}
\label{original}
For $0 \leq \theta < 1/2$, what is the least $f(\theta)$ such
that, for some $c_1, c_2 > 0$, any interval $[x - c_1 x^{f(\theta)}, x
+ c_1 x^{f(\theta)}]$ contains an integer $n$ with $n = ab$, and
$a$, $b$ are integers in the interval $[x^{1/2} - c_2 x^\theta,
x^{1/2} + c_2 x^\theta]$? Note: $c_1$ and $c_2$ may depend on
$\theta$.
\end{question}
We call such an integer $n$ an ``almost square" as it can be written as a
product of two integers about $\sqrt{n}$, with small error
$O(n^\theta)$. Clearly, $f$ is a non-increasing function and $0
\leq f(\theta) \leq 1/2$. Furthermore, we have
\begin{thm}
\label{thm1} For $0 \leq \theta < 1/4$, $f(\theta) = 1/2$.
\end{thm}
\begin{thm}
\label{thm2} For $0 \leq \theta < 1/2$, $f(\theta) \geq 1/2 -
\theta$.
\end{thm}
In fact, one suspects the following
\begin{conj}
\label{conj1} For $1/4 \leq \theta < 1/2$, $f(\theta) = 1/2 -
\theta$.
\end{conj}
Towards Conjecture \ref{conj1}, we have
\begin{thm}
\label{thm3} $f(1/4) = 1/4$.
\end{thm}
Conditionally, we have
\begin{thm}
\label{theorem3/10}
Assume Conjecture \ref{c2} on a certain average of twisted incomplete Sali\'{e} sum. For any $\epsilon > 0$, $f(\theta) \leq 1/2 - \theta + \epsilon$ for $1/4 < \theta \leq 3/10$.
\end{thm}

{\bf Notations:} $\{ x \} := x - [x]$ is the fractional part of
$x$. $\| x \| := \min_{l \in \mathbb{Z}} |x-l|$ is the distance
from $x$ to the closest integer. $f(x) = O(g(x))$ or $f(x) \ll g(x)$ mean
$|f(x)| \leq C g(x)$ for some constant $C > 0$. $f(x) \asymp g(x)$
means $g(x) \ll f(x) \ll g(x)$.

\bigskip

Our study is based on two ideas. First,
\begin{equation}
\label{key} ab = \Bigl(\frac{b+a}{2} - \frac{b-a}{2}\Bigr)
\Bigl(\frac{b+a}{2} + \frac{b-a}{2}\Bigr) =
\Bigl(\frac{b+a}{2}\Bigr)^2 - \Bigl(\frac{b-a}{2}\Bigr)^2
\end{equation}
which transforms the problem on factoring $n = ab$ into representing $n$ as the difference of two squares. Hence, for any $x$, we wish to find $ab = (\frac{b+a}{2})^2 - (\frac{b-a}{2})^2$ close to $x$ with $a, b$ close to $\sqrt{x}$. In other words,
we want $x + (\frac{b-a}{2})^2$ close to $(\frac{b+a}{2})^2$, or
equivalently, $\sqrt{x + (\frac{b-a}{2})^2}$ close to
$\frac{b+a}{2}$. Now suppose that $b-a$ is an even integer and
$d = \frac{b-a}{2}$. Thus, we transform Question \ref{original} to
\begin{question}
\label{transform}
Find integer $0 \leq d \leq c_2 x^\theta$ such that $\sqrt{x + d^2}$ is close to
an integer.
\end{question}
The second idea is the use of Taylor's expansion:
\begin{equation}
\label{taylor}
\sqrt{x + d^2} = \sqrt{x} \sqrt{1+\frac{d^2}{x}} = \sqrt{x} \Bigl[1 + \frac{d^2}{2x} + O\Bigl(\frac{d^4}{x^2}\Bigr) \Bigr] = \sqrt{x} + \frac{d^2}{2 \sqrt{x}} + O\Bigl(\frac{d^4}{x^{3/2}}\Bigr).
\end{equation}
This draws connection to the famous problem on the distribution of $n^2 \alpha \pmod 1$. In our situation, $\alpha = \frac{1}{2\sqrt{x}}$ and $0 \leq n \leq c_2 x^\theta$. The reader will see the reason for using Conjecture \ref{c2} in section \ref{n^2}. It may be worth mentioning here that the $3/10$ in Theorem \ref{theorem3/10} is due to the error term in (\ref{taylor}) and, ignoring this error term, Friedlander and Iwaniec's method in [\ref{FI}] only works for $\theta$ up to $1/3$.

\smallskip

Furthermore, from (\ref{key}), one may wonder if our problem is related to gaps between sums of two squares. In fact, using the same argument on $\sqrt{x - d^2}$ instead of $\sqrt{x + d^2}$, we have
\begin{thm}
Assume Conjecture \ref{c2}. For any $\epsilon > 0$, there exists some constant $c > 0$ such that, for any $x \geq 1$, the interval $[x, x + c x^{1/5+\epsilon}]$ always contains an integer which is the sum of two squares.
\end{thm}
In general, the same method gives similar result on gaps between values represented by a binary quadratic form. This can be done simply by completing the squares or using (\ref{key}) to transform a binary quadratic form into the form $AX^2 + BY^2$ (Note: One may need to impose certain restrictions on the ranges of the variables when the quadratic form is indefinite).

\smallskip

{\bf Acknowledgement} The author would like to thank Matthew Young for many stimulating discussions, and the American Institute of Mathematics for financial support.
\section{$0 \leq \theta < 1/4$ and Lower bound of $f$}
Proof of Theorem \ref{thm1}: Let $0 \leq \theta < 1/4$ and $c >
0$. It suffices to show that there exists arbitrarily large $x$
such that $|x - ab| \gg x^{1/2}$ for all integers $x^{1/2} -
cx^\theta \leq a \leq b \leq x^{1/2} + cx^\theta$. Now, consider any large $x$ with $\{\sqrt{x}\} = 1/4$. By (\ref{key}) and triangle inequality,
\begin{equation*}
\begin{split}
|x - ab| \geq& \Big|x - \Bigl(\frac{b+a}{2}\Bigr)^2\Big| -
\Big|\Bigl(\frac{b-a}{2}\Bigr)^2\Big| \\
\geq& \Big|\sqrt{x} - \frac{b+a}{2}\Big| \Big|\sqrt{x} +
\frac{b+a}{2}\Big| - c^2 x^{2\theta} \\
\geq& \frac{1}{4} \sqrt{x} - c^2 x^{2 \theta} \gg x^{1/2}.
\end{split}
\end{equation*}
Therefore, $f(\theta) \geq 1/2$ for $0 \leq \theta < 1/4$ which
gives the theorem as $f(\theta) \leq 1/2$ from consideration of 
perfect squares.

\bigskip

Proof of Theorem \ref{thm2}: For large $x$ and $x^{1/2} - c
x^\theta \leq a, b \leq x^{1/2} + c x^\theta$, there are at most
$4c^2 x^{2 \theta}$ distinct integers of the form $ab$ in the
interval $[x - c x^{1/2 + \theta}, x]$. Hence, there is a gap of
size at least 
$$\frac{c x^{1/2 + \theta}}{4 c^2 x^{2 \theta}} =
\frac{x^{1/2 - \theta}}{4c}$$
between two consecutive integers of the form $ab$. Now, pick $y$ as
the mid-point of this gap, then $[y - \frac{1}{10c} y^{1/2 - \theta}, y + \frac{1}{10c} y^{1/2 - \theta}]$ is an interval containing no integer of
the form $ab$ where $y^{1/2} - \frac{c}{2} y^\theta \leq a, b \leq y^{1/2} +
\frac{c}{2} y^\theta$. Since $c$ is arbitrary and $y$ can be
arbitrarily large, we must have $f(\theta) \geq 1/2 - \theta$.
\section{$\theta = 1/4$ and Conjecture \ref{conj1}}

From now on, we shall focus on $1/4 \leq \theta < 1/2$ as the case
$0 \leq \theta < 1/4$ is settled by Theorem \ref{thm1}.

Proof of Theorem \ref{thm3}: Consider the sequence ${\mathcal S} = \{\sqrt{x + d^2} \}_{cx^\theta < d \leq 2cx^\theta}$. The distance between successive
elements of ${\mathcal S}$ is
$$\sqrt{x + (d+1)^2} - \sqrt{x + d^2} = \frac{2d+1}{\sqrt{x + (d+1)^2}
+ \sqrt{x + d^2}} \asymp \frac{d}{\sqrt{x}} \asymp \frac{c}{x^{1/2
- \theta}}.$$
Since there are $c x^\theta + O(1)$ of these $d$'s, the distance between the first and last elements of ${\mathcal S}$ is
$$\asymp c x^\theta \frac{c}{x^{1/2-\theta}} = c^2 x^{2\theta -
1/2} > 1$$
when $\theta \geq 1/4$ and $c>1$. Hence, for some $cx^\theta < d_0 \leq 2cx^\theta$,
$$\|\sqrt{x + d_0^2}\| \ll \frac{c}{x^{1/2 - \theta}} \leq \frac{c}{x^{1/4}}.$$
Thus, for some integer $D$,
\begin{equation*}
\begin{split}
|\sqrt{x + d_0^2} - D| &\ll \frac{c}{x^{1/4}} \\
|x + d_0^2 - D^2| &\ll c x^{1/4} \\
|x - (D - d_0)(D + d_0)| &\ll c x^{1/4}.
\end{split}
\end{equation*}
Note that $\sqrt{x + c^2 x^{1/2}} + O(1) \leq D \leq \sqrt{x +
4c^2 x^{1/2}} + O(1)$ which implies $D = \sqrt{x} + O(c)$.
Therefore, $a = D-d_0$ and $b = D+d_0$ leads to $f(\theta) \leq 1/4$ for $1/4 \leq \theta < 1/2$. Combining this with Theorem \ref{thm2}, we have $f(1/4) = 1/4$.

\bigskip

Reason for Conjecture \ref{conj1}: One speculates that the numbers $\{\sqrt{x + d^2}\}_{0 \leq d \leq c x^\theta}$ are uniformly distributed $\pmod 1$. Thus, one expects that, for some $0 \leq d_0 \leq c x^\theta$,
$$\|\sqrt{x + d_0^2}\| \ll \frac{1}{c x^\theta}.$$
Then, for some integer $D$,
\begin{equation*}
\begin{split}
|\sqrt{x + d_0^2} - D| &\ll \frac{1}{c x^\theta} \\
|x + d_0^2 - D^2| &\ll \frac{1}{c} x^{1/2 - \theta} \\
|x - (D - d_0)(D + d_0)| &\ll \frac{1}{c} x^{1/2 - \theta}.
\end{split}
\end{equation*}
Now, $D = \sqrt{x+d_0^2} + O(1) = \sqrt{x} + O(d_0^2/\sqrt{x}) = \sqrt{x} + O(x^{2\theta - 1/2}) = \sqrt{x} + O(x^\theta)$. Therefore, $a = D-d_0$ and $b = D+d_0$ leads to $f(\theta) \leq 1/2 - \theta$.
\section{Connection to $n^2 \alpha \pmod 1$}
\label{n^2}
Hardy and Littlewood [\ref{HL}] conjectured that, for any real $\alpha$, there exists some $1 \leq n \leq N$ such that
$$\| n^2 \alpha \| \leq \frac{C}{N}$$
for some absolute constant $C > 0$. H. Heilbronn [\ref{H}] proved that
$$\| n^2 \alpha \| \ll_\epsilon \frac{1}{N^{1/2 - \epsilon}}$$
for any $\epsilon > 0$. The current best unconditional bound is due to A. Zaharescu [\ref{Z}], who showed that
$$\| n^2 \alpha \| \ll_\epsilon \frac{1}{N^{4/7 - \epsilon}}.$$
By assuming a certain conjecture on twisted incomplete Kloosterman sum, Friedlander and Iwaniec [\ref{FI}] proved that
$$\| n^2 \alpha \| \ll_\epsilon \frac{1}{n^{2/3 - \epsilon}}$$
for infinitely many positive integers $n$; and they claimed that, by a similar argument, one can get
$$\| n^2 \alpha \| \ll_\epsilon \frac{1}{n^{1 - \epsilon}}$$
assuming
\begin{conj}
\label{c1}
Let $a, q \geq 2$ be integers with $(a,q) = 1$ and $q$ not a perfect square. Let $H \geq 1$, $K \geq 1$ be reals. Then, for any $\epsilon > 0$,
\begin{equation*}
\begin{split}
& \mathop{\sum_{1 \leq h \leq H}}_{(h,q) = 1} \sum_{0 \leq k < K}  \Bigl(\frac{h}{q}\Bigr) e\Bigl(\frac{a \bar{h} k^2}{q}\Bigr) \\
& \ll_\epsilon (H^{1/2} K^{1/2} + H^{3/4} + K + q^{-1/2} HK + q^{-1/2} K^2) q^\epsilon
\end{split}
\end{equation*}
where $(\frac{\cdot}{\cdot})$ is the Jacobi symbol, $e(u) := e^{2\pi i u}$ and $\bar{h}$ denotes the multiplicative inverse of $h$ modulo $q$.
\end{conj}
Note: Friedlander and Iwaniec's results are weaker in the sense that they lack regularity on the occurrence of $n$, but they are better and best possible in terms of the exponent.

\smallskip

However, our situation in (\ref{taylor}) is a little different in two aspects: 

(i) we can no longer look at small distance from integers only because of the presence of $\sqrt{x}$. Instead, we need the fractional part of $\frac{d^2}{2\sqrt{x}}$ inside some small interval anywhere in $[0,1) \pmod 1$.

(ii) Our $\alpha = \frac{1}{2\sqrt{x}}$ has dependence on the parameter $x$ (and hence $N = [x^\theta]$).

Nevertheless, we shall use Friedlander and Iwaniec's method. Due to (i) and (ii), we need to assume something more general, namely,
\begin{conj}
\label{c2}
Let $a, q \geq 2$ be integers with $(a,q) = 1$ and $q$ not a perfect square. Let $H \geq 1$, $K \geq 1$ and $\lambda$, $\mu$ be reals. Then, for any $\epsilon > 0$,
\begin{equation*}
\begin{split}
& \mathop{\sum_{1 \leq h \leq H}}_{(h,q) = 1} e(\lambda h) \sum_{0 \leq k < K}   e(\mu k) \Bigl(\frac{h}{q}\Bigr) e\Bigl(\frac{a \bar{h} k^2}{q}\Bigr) \\
& \ll_\epsilon (H^{1/2} K^{1/2} + H^{3/4} + K + q^{-1/2} HK + q^{-1/2} K^2) q^\epsilon .
\end{split}
\end{equation*}
\end{conj}
\section{Gauss sums}
This section parallels section 2 of [\ref{FI}]. Let
$$G(a,b;q) := \sum_{d \pmod q} e \Bigl(\frac{ad^2 + bd}{q}\Bigr).$$ 
Lemma \ref{lemma1} and \ref{lemma2} are very similar to those in [\ref{FI}], using properties of Gauss sums and Conjecture \ref{c2}. We shall omit their proofs.
\begin{lem}
\label{lemma1}
Let $q$ be odd, $(p,q) = 1$, $H, K \geq 1$ and $\lambda$, $\mu$ real. We have
\begin{equation}
\label{l1}
\begin{split}
& \sum_{0 \leq k < K} e(\mu k) \mathop{\sum_{1 \leq h \leq H}}_{(h,q) = 1} e(\lambda h) G(hp, \pm k; q) \\
& \ll_\epsilon (q^{1/2} H^{1/2} K^{1/2} + q^{1/2} H^{3/4} + q^{1/2} K + HK + K^2) q^\epsilon .
\end{split}
\end{equation}
\end{lem}
\begin{lem}
\label{lemma2}
The estimate (\ref{l1}) of Lemma \ref{lemma1} still holds if the condition $(h,q) = 1$ is removed.
\end{lem}
Define, for $0 < \Delta < 1/2$ and $h \neq 0$,
\begin{equation}
\label{c}
c(h) := \Delta \Bigl(\frac{\sin \pi \Delta h}{\pi \Delta h}\Bigr)^2.
\end{equation}
\begin{lem}
\label{lemma3}
Let $q$ be odd, $(p,q) = 1$, $K \geq 1$, $\lambda$, $\mu$ real and $0 < \Delta < 1/2$. We have
\begin{equation*}
\begin{split}
& \sum_{0 \leq k < K} e(\mu k) \sum_{h = 1}^{\infty} c(h) e(\lambda h) G(hp, \pm k; q) \\
&\ll_\epsilon q^\epsilon \log{\frac{q}{\Delta}} (q^{1/2} K^{1/2} \Delta^{1/2} + q^{1/2} \Delta^{1/4} + q^{1/2} K \Delta + K + K^2 \Delta).
\end{split}
\end{equation*}
\end{lem}

Proof: Let $L \geq 1/\Delta$ be a parameter to be chosen later. We can split the above sum into
\begin{equation*}
\begin{split}
\sum_{0 \leq k < K} & e(\mu k) \sum_{1 \leq h \leq L} c(h) e(\lambda h) G(hp, \pm k; q) \\
+& \sum_{0 \leq k < K}  e(\mu k) \sum_{h > L} c(h) e(\lambda h) G(hp, \pm k; q)
= \Sigma_1 + \Sigma_2.
\end{split}
\end{equation*}
Like [\ref{FI}], we apply the trivial estimates
\begin{equation}
\label{trivial}
c(h) \ll \Delta^{-1} h^{-2} \; \mbox{ and } \; G(hp, k; q) \ll q
\end{equation}
to $\Sigma_2$ and get
\begin{equation}
\label{3.1}
\Sigma_2 \ll K \int_{L}^{\infty} \frac{1}{\Delta h^2} q \; dh \ll \frac{Kq}{\Delta L}.
\end{equation}
For $\Sigma_1$, we apply Lemma \ref{lemma2} and partial summation getting
\begin{equation}
\label{3.2}
\begin{split}
\Sigma_1 =& \int_{1}^{L} c(H) d\Bigl(\sum_{1 \leq h \leq H} e(\lambda h) \sum_{0 \leq k < K} e(\mu k) G(hp, \pm k; q) \Bigr) \\
\ll_\epsilon& c(L) (q^{1/2} L^{1/2} K^{1/2} + q^{1/2} L^{3/4} + q^{1/2} K + LK + K^2) q^\epsilon \\
&+ q^\epsilon \int_{1}^{L} |c'(H)| (q^{1/2} H^{1/2} K^{1/2} + q^{1/2} H^{3/4} + q^{1/2} K + HK + K^2) dH \\
\ll_\epsilon& \frac{q^\epsilon}{\Delta L} \Bigl(\frac{q^{1/2} K^{1/2}}{L^{1/2}} + \frac{q^{1/2}}{L^{1/4}} + \frac{q^{1/2} K}{L} + K + \frac{K^2}{L} \Bigr) \\
&+ q^\epsilon \int_{1}^{1/\Delta} \frac{\Delta}{H} (q^{1/2} H^{1/2} K^{1/2} + q^{1/2} H^{3/4} + q^{1/2} K + HK + K^2) dH \\
&+ q^\epsilon \int_{1/\Delta}^{L} \frac{1}{H^2} (q^{1/2} H^{1/2} K^{1/2} + q^{1/2} H^{3/4} + q^{1/2} K + HK + K^2) dH
\end{split}
\end{equation}
as $|c'(H)| \ll H^{-2} |\sin{(\pi \Delta H)}|$. By simple integration and estimation, one has the above two integrals
\begin{equation}
\label{3.3}
\ll_\epsilon q^\epsilon \log L (q^{1/2} K^{1/2} \Delta^{1/2} + q^{1/2} \Delta^{1/4} + q^{1/2} K \Delta + K + K^2 \Delta)
\end{equation}
Now, we pick $L = q/\Delta$ and the lemma follows from (\ref{3.1}), (\ref{3.2}) and (\ref{3.3}).
\begin{lem}
\label{lemma4}
Let $\{ a_k \}_{k=0}^{\infty}$ be a sequence of complex numbers and $f(x)$ be some positive function. If $\sum_{k < K} a_k \ll f(K)$ for all $K > 0$, then
$$\sum_{K \leq k < 2K} \frac{a_k}{k^2} \ll \frac{f(2K)}{K^2} + \int_{K}^{2K} \frac{f(u)}{u^3} du.$$
\end{lem}

Proof: By partial summation,
\begin{equation*}
\sum_{K \leq k < 2K} \frac{a_k}{k^2} = \int_{K}^{2K} \frac{1}{u^2} d \sum_{k < u} a_k
\ll \frac{f(2K)}{(2K)^2} + 2\int_{K}^{2K} \frac{\sum_{k < u} a_k}{u^3} du
\end{equation*}
which gives the lemma.
\begin{lem}
\label{lemma5}
Let $\{ a_k \}_{k=0}^{\infty}$ be a sequence of complex numbers and $f(x)$ be some positive function. If $\sum_{k < K} a_k e(\delta k) \ll f(K)$ for all $K > 0$ and all real $\delta$, then $\sum_{k < K} a_k \sin^2 \pi \delta k \ll f(K)$.
\end{lem}

Proof: Use $\sin^2 \pi \delta \theta = \frac{1}{2} (1 - \cos 2\pi \delta \theta)$, $\cos 2\pi \delta \theta = \frac{1}{2} (e(\delta \theta) + e(-\delta \theta))$ and triangle inequality.
\section{Towards the proof of Theorem \ref{theorem3/10}}
Let $q$ be odd and not a perfect square, $(p, q)=1$, $0 < \Delta < 1/2$ and $\lambda$ be any real number. Consider
\begin{equation}
\label{Slambda}
S_{p/q}^\lambda (\Delta, N) := \# \Bigl\{1 \leq n \leq N : \{ \frac{p}{q} n^2 \} \in (\lambda - \Delta, \lambda + \Delta) \pmod 1 \Bigr\}.
\end{equation}
\begin{lem}
\label{keylemma}
For $q > 2N$ and $N \geq 1$, we have
$$S_{p/q}^\lambda (\Delta, N) \geq \frac{1}{2} \Delta N + O\Bigl(q^\epsilon \log^2 \frac{q}{\Delta} \Bigl[N^{1/2} \Delta^{1/2} + \frac{N \Delta^{1/4}}{q^{1/2}} + \Delta q^{1/2} + 1 + \frac{q \Delta}{N} + \frac{N^5}{\Delta q^4}\Bigr] \Bigr).$$
\end{lem}

Proof: Let
\[ f(x) := \Bigl\{ \begin{array}{ll}
x + 1, & \mbox{ if } -1 \leq x \leq 0, \\
1 - x, & \mbox{ if } 0 \leq x \leq 1, \\
0, & \mbox{ otherwise;}
\end{array} \]
\[ t(x) := \Bigl\{ \begin{array}{ll}
x/\Delta + 1, & \mbox{ if } -\Delta \leq x \leq 0, \\
1 - x/\Delta, & \mbox{ if } 0 \leq x \leq \Delta, \\
0, & \mbox{ otherwise;}
\end{array} \]
\[ g(x) := \sum_{n = -\infty}^{\infty} t(x-n); \]
and
\[ g_\lambda(x) := g(x - \lambda). \]
Note: $f(x)$ and $g(x)$ are just the same as those in [\ref{FI}] and $g_\lambda(x)$ is a shift of $g(x)$ to the right by $\lambda$. Then, one has
\begin{equation}
\label{count}
1 + 2S_{p/q}^\lambda (\Delta, N) \geq \sum_{n = -\infty}^{\infty} f\Bigl(\frac{n}{N}\Bigr) g_\lambda \Bigl(\frac{p}{q} n^2\Bigr).
\end{equation}
The function $g_\lambda$ has Fourier expansion
$$g_\lambda(x) = \sum_{h = -\infty}^{\infty} c(h) e(-\lambda h) e(h x)$$
where $c(h)$ is defined by (\ref{c}). Using this, the right hand side of (\ref{count})
\begin{equation*}
\begin{split}
=& \sum_{h = -\infty}^{\infty} c(h) e(-\lambda h) \sum_{n = -\infty}^{\infty} e\Bigl(\frac{hp n^2}{q}\Bigr) f\Bigl(\frac{n}{N}\Bigr) \\
=& c(0) \sum_{n = -\infty}^{\infty} f\Bigl(\frac{n}{N}\Bigr) + \sum_{h \neq 0} c(h) e(-\lambda h) \sum_{n = -\infty}^{\infty} e\Bigl(\frac{hp n^2}{q}\Bigr) f\Bigl(\frac{n}{N}\Bigr) = \Delta N + R
\end{split}
\end{equation*}
for $N$ an integer and $R$ being the sum over $h \neq 0$. Now, we just follow the same calculation of applying Poisson summation in [\ref{FI}]. The inner sum of $R$ becomes
\begin{equation*}
\begin{split}
\sum_n =& \sum_{d\pmod q} e\Bigl(\frac{hp d^2}{q}\Bigr) \sum_{n \equiv d\pmod q} f\Bigl(\frac{n}{N}\Bigr) \\
=& \frac{N}{q} \sum_{k} G(hp, k; q) \hat{f}\Bigl(\frac{kN}{q}\Bigr)
\end{split}
\end{equation*}
where $\hat{f}(y) = (\frac{\sin{\pi y}}{\pi y})^2$ is the Fourier transform of $f$. Then
\begin{equation}
\begin{split}
\label{middle}
R =& \frac{N}{q} \sum_k \Bigl(\frac{\sin{(\pi k N/q)}}{\pi k N/q}\Bigr)^2 \sum_{h \neq 0} c(h) e(-\lambda h) G(hp, k; q) \\
=& \frac{N}{q} (\sum_{k < q/N} + \sum_{q/N \leq k < M} + \sum_{k \geq M})
\end{split}
\end{equation}
where $M = (q/N)^\sigma$ and $\sigma > 1$ is some parameter to be chosen later. We apply the trivial estimates in (\ref{trivial}) to bound
\begin{equation}
\label{s3}
\sum_{k \geq M} \ll \frac{q^2}{N^2} \sum_{k \geq M} \frac{1}{k^2} \sum_{h = 1}^{\infty} \frac{q}{\Delta h^2} \ll \frac{q^3}{\Delta N^2 M}.
\end{equation}
By Lemma \ref{lemma3} with $\mu = 0$; $K=1$ (to deal with $k=0$) and $K=q/N$, and partial summation, we have
\begin{equation}
\label{s1}
\begin{split}
\sum_{k < q/N} =& \int_{1}^{q/N} \Bigl(\frac{\sin{(\pi u N/q)}}{\pi u N/q}\Bigr)^2 d\Bigl( \sum_{0 \leq k < u} \sum_{h \neq 0} c(h) e(-\lambda h) G(hp, k; q) \Bigr) \\
&+ O\Bigl( q^\epsilon \log \frac{q}{\Delta} (q^{1/2} \Delta^{1/4} + 1) \Bigr) \\
\ll_\epsilon& q^\epsilon \log \frac{q}{\Delta} (q^{1/2} (\frac{q}{N})^{1/2} \Delta^{1/2} + q^{1/2} \Delta^{1/4} + q^{1/2} \frac{q}{N} \Delta + \frac{q}{N} + (\frac{q}{N})^2 \Delta) \\
&+ 2 \int_{1}^{q/N} \frac{|\sin{(\pi u N/q)}|}{N u^2/q} \Big|\sum_{0 \leq k < u} \sum_{h=1}^{\infty} c(h) e(-\lambda h) G(hp, k; q) \Big| du \\
&+ q^\epsilon \log \frac{q}{\Delta} (q^{1/2} \Delta^{1/4} + 1) \\
\ll_\epsilon& q^\epsilon \log \frac{q}{\Delta} \Bigl(\frac{q \Delta^{1/2}}{N^{1/2}} + q^{1/2} \Delta^{1/4} \log \frac{q}{N} + \frac{q^{3/2} \Delta}{N} + \frac{q}{N} + \frac{q^2 \Delta}{N^2}\Bigr)
\end{split}
\end{equation}
Suppose $K$ is a positive integer satisfying $2^K \frac{q}{N} = (\frac{q}{N})^\sigma = M$ (hence we shall choose $K$ at the end).
Let
$$a_k := \sum_{h \neq 0} c(h) e(-\lambda h) G(hp, k; q),$$
$$B(x) := q^\epsilon \log{\frac{q}{\Delta}} \Bigl(q^{1/2} x^{1/2} \Delta^{1/2} + q^{1/2} \Delta^{1/4} + q^{1/2} x \Delta + x + x^2 \Delta\Bigr),$$
and $\delta := N/q$. Applying Lemma \ref{lemma3} to Lemma \ref{lemma5} and then Lemma \ref{lemma4} with these $a_k$ and $\delta$, we have
\begin{equation}
\label{s2}
\begin{split}
\sum_{q/N \leq k < M} &= \sum_{i=1}^{K} \sum_{2^{i-1} q/N \leq k < 2^i q/N} \Bigl(\frac{\sin{(\pi u N/q)}}{\pi u N/q}\Bigr)^2 a_k \\
\ll& \sum_{i=1}^{K} \Bigl(\frac{q}{N}\Bigr)^2 \Bigl[\frac{B(2^i q/N)}{(2^i q/N)^2} + \int_{2^{i-1}q/N}^{2^i q/N} \frac{B(x)}{x^3} dx \Bigr] \\
\ll_\epsilon& \sum_{i=1}^{K} \Bigl(\frac{q}{N}\Bigr)^2 q^\epsilon \log \frac{q}{\Delta} \Bigl[\frac{q^{1/2} \Delta^{1/2}}{(2^i q/N)^{3/2}} + \frac{q^{1/2} \Delta^{1/4}}{(2^i q/N)^2} + \frac{q^{1/2} \Delta}{2^i q/N} + \frac{1}{2^i q/N} + \Delta \Bigr] \\
\ll_\epsilon& q^\epsilon \log \frac{q}{\Delta} \Bigl[\frac{q \Delta^{1/2}}{N^{1/2}} + q^{1/2} \Delta^{1/4} + \frac{\Delta q^{3/2}}{N} + \frac{q}{N} + \frac{q^2 \Delta}{N^2} (K + \log \frac{q}{N}) \Bigr].
\end{split}
\end{equation}
Putting (\ref{s3}), (\ref{s1}) and (\ref{s2}) into (\ref{middle}), we have
$$R \ll_\epsilon q^\epsilon \log \frac{q}{\Delta} (K + \log \frac{q}{N}) \Bigl[N^{1/2} \Delta^{1/2} + \frac{N \Delta^{1/4}}{q^{1/2}} + \Delta q^{1/2} + 1 + \frac{q \Delta}{N} + \frac{q^2}{\Delta N M} \Bigr]$$
which gives the lemma by choosing $K = [\frac{5 \log{(q/N)}}{\log 2}]$.
\section{Completion of the proof of Theorem \ref{theorem3/10}}
Let $\lambda$ be any real number. Based on (\ref{taylor}), we want to get a positive lower bound for
$$S_{1/2\sqrt{x}}^{\lambda}(\Delta, N) := \# \Bigl\{1 \leq n \leq N : \{ \frac{1}{2\sqrt{x}} n^2 \} \in (\lambda - \Delta, \lambda + \Delta) \pmod 1 \Bigr\}$$
(defined similar to (\ref{Slambda})) for appropriate choices of $\Delta$ and $N$. We shall pick $N = [x^\theta]$ and $\Delta = 1/x^{\theta - 2\epsilon}$ for small $\epsilon > 0$. In the notation of (\ref{Slambda}), we pick $p = 1$ and $q =$ $2[\sqrt{x}] + 1$ or $2[\sqrt{x}] + 3$ so that $q$ is not a square. Clearly, $(p,q) = 1$. Thus,
$$\Big|\frac{1}{2\sqrt{x}} - \frac{1}{q}\Big| \leq \frac{1}{x},$$
$$\Big|\frac{1}{2\sqrt{x}} n^2 - \frac{n^2}{q}\Big| \leq \frac{n^2}{x} \leq \frac{1}{x^{1-2\theta}}.$$
Therefore,
$$S_{1/2\sqrt{x}}^\lambda (\Delta, N) \geq S_{1/q}^\lambda (\frac{\Delta}{2}, N)$$
provided $1/x^{1-2\theta} \leq \Delta/2$. This is satisfied by our choice of $\Delta$ above when $1/4 < \theta < 1/3$ (with $\epsilon$ small enough). Now, we apply Lemma \ref{keylemma}. With our choices of $\Delta$, $N$ and $q$, one can check that
$$ q^{2\epsilon} \Bigl[N^{1/2} \Delta^{1/2} + \frac{N \Delta^{1/4}}{q^{1/2}} + \Delta q^{1/2} + 1 + \frac{q \Delta}{N} + \frac{N^5}{\Delta q^4}\Bigr] \ll \Delta N.$$
Hence, for $x$ large enough,
$$S_{1/2\sqrt{x}}^\lambda (\Delta, N) \geq \frac{1}{8} \Delta N > 0.$$
Consequently, we can conclude that, for any $\lambda$, there is some $0 \leq d \leq N$ such that the fractional part, $\{ \frac{d^2}{2\sqrt{x}} \}$, lies in $(\lambda - \Delta, \lambda + \Delta) \pmod 1$. In particular, when $1/4 < \theta < 1/3$, we can find $0 \leq d \leq x^\theta$ such that
$$\Big\| \sqrt{x} + \frac{d^2}{2\sqrt{x}} \Big\| \ll_{\epsilon, \theta} \frac{1}{x^{\theta - \epsilon}}$$
for any $\epsilon > 0$. By (\ref{taylor}), this implies
$$\Big\| \sqrt{x + d^2} \Big\| \ll \frac{1}{x^{\theta - \epsilon}} + \frac{x^{4\theta}}{x^{3/2}}.$$
When $1/4 < \theta \leq 3/10$, the second error term is smaller than the first error term above. Thus, for $1/4 < \theta \leq 3/10$, there exist integers $0 \leq d \leq x^\theta$ and $D$ such that
\begin{equation*}
\begin{split}
& |\sqrt{x + d^2} - D| \ll_{\epsilon, \theta} \frac{1}{x^{\theta - \epsilon}}, \\
& |x + d^2 - D^2| \ll_{\epsilon, \theta} x^{1/2 - \theta + \epsilon} \\
& |x - (D-d)(D+d)| \ll_{\epsilon, \theta} x^{1/2 - \theta + \epsilon}.
\end{split}
\end{equation*}
Note: $D = \sqrt{x+d^2} + O(1) = \sqrt{x} + O(d^2/\sqrt{x}) = \sqrt{x} + O(x^{2\theta - 1/2}) = \sqrt{x} + O(x^{\theta})$. Consequently, $(D-d)(D+d)$ is the required integer for our theorem.
\section{$\theta = 1/2$}

So far, we have excluded the endpoint $\theta = 1/2$. One reason is that Conjecture \ref{conj1} cannot be extended to $\theta = 1/2$. The situation is more delicate. We have
\begin{thm}
Let $F(x)$ be a positive real value function such that any interval $[x - F(x), x + F(x)]$ contains an integer $n$ with $n = ab$; $a$, $b$ are integers in the interval $[0, c \sqrt{x}]$ for some constant $c > 0$. Then
$$(\log x)^{0.086} \ll F(x) \ll x^{1/4}.$$
\end{thm}

Proof: The upper bound follows from Theorem \ref{thm3}. (Note: Assuming Conjecture \ref{c2}, the upper bound can be improved to $x^{1/5 + \epsilon}$ by Theorem \ref{theorem3/10}.) As for the lower bound, let us consider the interval $I = [x, 2x]$. By the definition of $F$, there are $\gg x/F(x)$ distinct integers of the form $ab$ with $0 \leq a, b \leq c \sqrt{2x}$ in $I$. But, by Erd\"{o}s [\ref{E}], there are $\ll c^2 2x / (\log x)^{\alpha +o(1)}$ distinct integers from the ``multiplication table" of size $c \sqrt{2x} \times c \sqrt{2x}$ where $\alpha = 1 - \log{(e \log 2)}/\log 2 = 0.0860...$ . Hence
$$\frac{x}{F(x)} \ll \frac{c^2 x}{(\log x)^{0.086}}$$
which gives the lower bound.

The situation here is similar to that on gaps between sums of two squares. Both best upper bounds, so far, are obtained by elementary means. In fact, one may guess $F(x) \ll x^\epsilon$ for every $\epsilon > 0$ (true under Conjecture \ref{conj1}), and even perhaps $F(x) \ll (\log x)^C$ for some $C > 0$.
\bigskip


Tsz Ho Chan\\
American Institute of Mathematics\\
360 Portage Avenue\\
Palo Alto, CA 94306\\
USA\\
thchan@aimath.org

\end{document}